\documentclass[reqno, 10pt]{amsart}
\usepackage{todonotes}
\presetkeys{todonotes}{inline}{}
\makeatletter
\providecommand\@dotsep{5}
\makeatother
\usepackage{amssymb, amsmath, amsthm, enumerate,  marginnote, mathabx, mathrsfs, mathtools, tikz, bm}
\usetikzlibrary{arrows,calc}
\usetikzlibrary{arrows.meta, cd, patterns}
\usepackage{bbm} 
\usepackage{esint} 
\usepackage{verbatim}

\usepackage{hyperref}
\hypersetup{
    colorlinks=true,
    linkcolor=blue,
    filecolor=magenta,
}

\usepackage[normalem]{ulem}

\usetikzlibrary{decorations.pathreplacing}
\usetikzlibrary{decorations,calligraphy}

\usepackage{todonotes}

\usepackage{arydshln, multirow} 

\usepackage{comment}
\theoremstyle{plain}
\numberwithin{equation}{section}
\newtheorem{theorem}{Theorem}[section]

\newtheorem{lemma}[theorem]{Lemma}
\newtheorem{proposition}[theorem]{Proposition}
\newtheorem{definition}[theorem]{Definition}

\newtheorem{conjecture}[theorem]{Conjecture}

\theoremstyle{definition}

\newtheorem*{induction hypothesis}{Induction hypothesis}

\theoremstyle{remark}

\def\inn#1#2{\langle#1, #2\rangle}

\newcommand{\supp}{\mathrm{supp}\,}

\newcommand{\R}{\mathbb{R}}
\newcommand{\Z}{\mathbb{Z}}
\newcommand{\N}{\mathbb{N}}

\newcommand{\ud}{\mathrm{d}}

\def\lc{\lesssim}

\newcommand{\Be}{\begin{equation}}
\newcommand{\Ee}{\end{equation}}

\newcommand{\Bm}{\begin{multline}}
\newcommand{\Em}{\end{multline}}
\newcommand{\Bea}{\begin{eqnarray}}
\newcommand{\Eea}{\end{eqnarray}}
\newcommand{\Beas}{\begin{eqnarray*}}
\newcommand{\Eeas}{\end{eqnarray*}}
\newcommand{\Benu}{\begin{enumerate}}
\newcommand{\Eenu}{\end{enumerate}}
\newcommand{\Bi}{\begin{itemize}}
\newcommand{\Ei}{\end{itemize}}

\def\intslash{\rlap{\kern  .32em $\mspace {.5mu}\backslash$ }\int}
\def\qsl{{\rlap{\kern  .32em $\mspace {.5mu}\backslash$ }\int_{Q_x}}}

\def\N{\mathbb N}

\def\emph#1{{\it #1 }}

\def\inn#1#2{\langle#1,#2\rangle}

\def\lc{\lesssim}

\newcommand{\la}{\lambda}

\def\fN{{\mathfrak {N}}}

\def\bbN{{\mathbb {N}}}

\def\bbR{{\mathbb {R}}}

\def\cS{{\mathcal {S}}}

\def\cZ{{\mathcal {Z}}}

\title[Local smoothing conjecture for curves]{A counterexample for local smoothing for averages over curves}

\author[D. Beltran]{David Beltran}
\author[J. Hickman]{Jonathan Hickman}
\date{\today}

\address{David Beltran: Departament d’An\`alisi Matem\`atica, Universitat de Val\`encia, Dr. Moliner 50, 46100 Burjassot, Spain}
\email{david.beltran@uv.es}
\address{Jonathan Hickman: School of Mathematics, James Clerk Maxwell Building, The King's Buildings, Peter Guthrie Tait Road, Edinburgh, EH9 3FD, UK.}
\email{jonathan.hickman@ed.ac.uk}


\begin{document}

\begin{abstract}
We provide a new necessary condition for local smoothing estimates for the averaging operator defined by convolution with a measure supported on a smooth non-degenerate curve in $\R^n$ for $n \geq 3$. This demonstrates a limitation in the strength of local smoothing estimates towards establishing bounds for the corresponding maximal functions when $n \geq 5$.
\end{abstract}

\maketitle




\section{Introduction} 




\subsection{Main result} For $n\ge 2$ let $\gamma \colon I \to \R^n$ be a smooth curve, 
where $I \subset \R$ is a compact interval, and $\chi \in C^{\infty}(\R)$ be a bump function supported on the interior of $I$. Given $t>0$, define the measure $\mu_t$ supported on $t$-dilates of $\gamma$ by
\begin{equation}\label{eq: measure}
    \inn{\mu_t}{g} = \int_\R g(t \gamma (s)) \chi(s) \, \ud s
\end{equation}
and consider the associated averaging operator
\begin{equation}\label{eq:averaging operator}
    A_tf(x) := \mu_t \ast f (x) = \int_{\R} f(x - t \gamma(s))\,\chi(s) \, \ud s.
\end{equation}
We will focus on these averaging operators for \textit{non-degenerate} curves: that is, smooth curves $\gamma: I \to \R^n$ for which there is a constant $c_0 > 0$ such that 
\begin{equation}\label{eq:nondegenerate}
    |\det(\gamma'(s), \cdots, \gamma^{(n)}(s))| \geq c_0 \qquad \textrm{for all $s \in I$}.
\end{equation}

We begin by discussing $A_t$ for a fixed value of $1 \leq t \leq 2$. 
Under the nondegeneracy hypothesis \eqref{eq:nondegenerate}, Ko, Lee and Oh \cite{KLO_Sob} proved that if $2(n-1)< p < \infty$, then
\begin{equation}\label{eq:fixed-time}
   \|A_t f\|_{L^p_{1/p}(\R^n)} \lesssim_{p,\gamma, \chi} \|f\|_{L^p(\R^n)}.
\end{equation}
This result is essentially sharp. First, it is well-known that $1/p$ is the best possible order of smoothing in \eqref{eq:fixed-time}. Second, it was shown in \cite{BGHS_Sob} that \eqref{eq:fixed-time} fails for $p < 2(n-1)$. More precisely, if
$2 \leq p \leq \infty$,  the curve $\gamma \colon I \to \R^n$ is non-degenerate and the inequality 
\begin{equation}\label{eq:fixed-time nec}
    \| A_t f \|_{L^p_\alpha(\R^n)} \lesssim_{p,\gamma,\chi} \| f \|_{L^p(\R^n)}
\end{equation}
holds, then we must have $\alpha \leq \min \{ \frac{1}{n}(\frac{1}{2} + \frac{1}{p} ), \frac{1}{p} \}$. 
The necessary condition $ \alpha \leq \frac{1}{n}(\frac{1}{2}+\frac{1}{p})$ was first observed for the helix $\gamma(s)=(\cos s, \sin s, s)$ in \cite{OS1999}; this was later extended to higher dimensions in \cite{BGHS_Sob}, where the sharp examples were also contextualised in relation to decoupling theory.

Here we are interested in the extent to which it is possible to improve the inequality \eqref{eq:fixed-time} by integrating locally in the $t$-variable. Our main result is the following.

\begin{theorem}\label{thm:LS nec} Let $2 \leq p \leq \infty$. If $\gamma \colon I \to \R^n$ is non-degenerate and the inequality 
\begin{equation}\label{eq:LS thm}
    \Big( \int_1^2\| A_t f \|_{L^p_\sigma(\R^n)}^p \, \ud t \Big)^{1/p} \lesssim_{p,\gamma,\chi} \| f \|_{L^p(\R^n)}
\end{equation}
holds, then we must have $\sigma \leq \sigma(p,n):= \min \Big\{ \frac{1}{n}, \frac{1}{n}\Big(\frac{1}{2} + \frac{2}{p} \Big), \frac{2}{p} \Big\}$.
\end{theorem}

Note that the critical index $\sigma(p,n)$ can be expressed as 
\begin{equation*}
    \sigma(p,n)=
    \begin{cases}
        \frac{1}{n} \quad & \text{ if }\,\, 2 \leq p \leq 4 , \\
        \frac{1}{n}(\frac{1}{2}+\frac{2}{p}) \quad & \text{ if }\,\, 4 \leq p \leq 4(n-1) , \\
        \frac{2}{p} \quad & \text{ if }\,\, 4(n-1) \leq p \leq \infty.
    \end{cases}
\end{equation*}
The condition $\sigma \leq \min \{\frac{1}{n}, \frac{2}{p} \}$ was already observed in \cite{KLO_Sob}. Moreover, on \cite[p.3]{KLO_Sob} the authors remark that `it seems to be plausible to conjecture' that \eqref{eq:LS thm} holds for $\sigma < \min \{\frac{1}{n}, \frac{2}{p} \}$ if $2 \leq p \leq \infty$. Thus, our main contribution is the additional condition $\sigma \leq \frac{1}{n}(\frac{1}{2} + \frac{2}{p} )$, which provides a counterexample to the above conjecture. We remark that our condition is only relevant for $n \geq 3$.\medskip

Theorem~\ref{thm:LS nec} provides a local smoothing variant of the necessary condition $\alpha \leq \frac{1}{n}(\frac{1}{2}+\frac{1}{p})$ for the fixed-time $L^p$-Sobolev estimate \eqref{eq:fixed-time nec}. Moreover, the proof is a direct modification of the construction used in \cite[Proposition 3.3]{BGHS_Sob}.




\subsection{A revised conjecture and consequences for the theory of geometric maximal operators} In light of Theorem~\ref{thm:LS nec}, one may be tempted to refine the conjectural bounds as follows.

\begin{conjecture}[Local smoothing conjecture for curve averages]\label{conj:LS}
    Let $n \geq 2$ and $2 < p < \infty$. If $\gamma:I \to \R^n$ is a non-degenerate curve, then the inequality \eqref{eq:LS thm} holds for all $\sigma < \sigma(p,n)$.
\end{conjecture}

This conjecture is solved affirmatively for $n=2$ combining the works \cite{GWZ, MSS1992}: in this context, it is essentially equivalent to the local smoothing problem for the wave equation in $\R^2$ posed by Sogge \cite{Sogge1991}, from which we borrow the terminology \textit{local smoothing}. For $n \geq 3$, Ko, Lee and Oh \cite{KLO_Sob} have verified Conjecture \ref{conj:LS} if $p > 4n-2$. The range $2 < p \leq 4n-2$ remains open.\medskip 

A major motivation for the study of local smoothing estimates of the type \eqref{eq:LS thm} is that they typically imply $L^p$-bounds for the maximal function 
\begin{equation}\label{eq: max fn}
    M_\gamma f(x):=\sup_{t>0} |A_t f(x)|. 
\end{equation}
This is a prototypical example of a geometric maximal operator, providing a natural generalisation of Bourgain's circular maximal function \cite{Bourgain86} to higher dimensions. 

\begin{conjecture}\label{conj:max funct}
Let $n\geq 2$ and $\gamma:I \to \R^n$ a non-degenerate curve. Then $M_\gamma$ maps $L^p(\R^n) \to L^p(\R^n)$ boundedly if and only if $p>n$.   
\end{conjecture}

This conjecture is known for $n=2$ \cite{Bourgain86, MSS1992} and $n=3$ \cite{PS2007, BGHS_Hel, KLO_Hel}. For $n \geq 4$ it is currently open, but it was shown in \cite{KLO_Sob} that $L^p(\R^n)$ bounds hold for $p>2(n-1)$.\medskip

Bounds for the maximal function \eqref{eq: max fn} can be proved using local smoothing estimates \eqref{eq:LS thm} via a standard argument involving Sobolev embedding and Littlewood--Paley theory (see, for instance, \cite[Proposition 2.1]{BGHS_Hel}). In particular, if \eqref{eq:LS thm} holds for $\sigma > 1/p$, then $M_\gamma$ is bounded on $L^p(\R^n)$. Many partial results towards Conjecture~\ref{conj:max funct} follow this proof strategy: local smoothing estimates are proved, which are then translated into maximal estimates via Sobolev embedding. This is the case for the current best bounds for $n \geq 4$ from \cite{KLO_Sob}.\footnote{We remark that Bourgain's argument \cite{Bourgain86} can also be reinterpreted in this paradigm: see \cite{Sogge1991}. The method used to study the $n=3$ case in \cite{KLO_Hel} relies on linearising the maximal function and proving $L^p$ estimates for $A_tf$ with respect to fractal measures; however, much of the technology used to prove the fractal estimates is at the level of local smoothing.} More precisely, the maximal estimates in \cite{KLO_Sob} were obtained as a consequence of the aforementioned local smoothing estimates \eqref{eq:LS thm} for $\sigma < 2/p$ and $p>4n-2$: interpolation of the local smooth estimates with a trivial $L^2 \to L^2_{1/n}$-bound gives then \eqref{eq:LS thm} with $\sigma > 1/p$ for $p > 2(n-1)$. \medskip

An interesting consequence of Theorem \ref{thm:LS nec} is that for $n \geq 4$ the inequality \eqref{eq:LS thm} can only hold for $\sigma > 1/p$ if $p > 2(n-2)$. This means that for $n \geq 5$ one cannot verify Conjecture~\ref{conj:max funct} in the whole range $p > n$ from sharp local smoothing via the usual argument. We feel that this observation is significant: the local smoothing approach has dominated work towards Conjecture~\ref{conj:max funct}, so Theorem~\ref{thm:LS nec} highlights an inherent limitation in much of our current understanding of the problem.





\subsection*{Notational conventions} Given a (possibly empty) list of objects $L$, for real numbers $A_p, B_p \geq 0$ depending on some Lebesgue exponent $p$ or dimension parameter $n$ the notation $A_p \lesssim_L B_p$, $A_p = O_L(B_p)$ or $B_p \gtrsim_L A_p$ signifies that $A_p \leq CB_p$ for some constant $C = C_{L,p,n} \geq 0$ depending on the objects in the list, $p$ and $n$. In addition, $A_p \sim_L B_p$ is used to signify that both $A_p \lesssim_L B_p$ and $A_p \gtrsim_L B_p$ hold. 
The length of a multiindex $\alpha\in \bbN_0^n$ is given by $|\alpha|=\sum_{i=1}^n{\alpha_i}$.

\subsection*{Acknowledgements}

The first author is supported by the grants RYC2020-029151-I and PID2022-140977NA-I00 funded by MICIU/AEI/10.13039/501100011033, ``ESF Investing in your future" and FEDER, UE. The second author is supported by New Investigator Award UKRI097.




\section{Preliminaries}




\subsection{Reduction to perturbations of the moment curve}
We begin with some standard reductions which have appeared frequently in the literature. A prototypical example of a smooth curve satisfying the non-degeneracy condition \eqref{eq:nondegenerate} is the \textit{moment curve} $\gamma_{\circ} \colon \R \to \R^n$, given by
\begin{equation*}
    \gamma_{\circ}(s) := \Big(s, \frac{s^2}{2}, \dots, \frac{s^n}{n!} \Big). 
\end{equation*}

We consider a class of model curves which are perturbations of $\gamma_\circ$.

\begin{definition} Given $n \geq 2$ and $0 < \delta < 1$, let $\mathfrak{G}_n(\delta)$ denote the class of all smooth curves $\gamma \colon [-1, 1] \to \R^n$ that satisfy the following conditions: 
\begin{enumerate}[i)]
    \item $\gamma(0) = 0$ and $\gamma^{(j)}(0) = \vec{e}_j$ for $1 \leq j \leq n$;
    \item $\|\gamma - \gamma_{\circ}\|_{C^{n+1}([-1,1])} \leq \delta$.
\end{enumerate}
Here $\vec{e}_j$ denotes the $j$th standard Euclidean basis vector and
\begin{equation*}
    \|\gamma\|_{C^{n+1}(I)} := \max_{1 \leq j \leq n + 1} \sup_{s \in I} |\gamma^{(j)}(s)| \qquad \textrm{for all $\gamma \in C^{n+1}(I;\R^n)$.}
\end{equation*}
\end{definition}

By Taylor expansion and standard scaling arguments, one can reduce the problem of studying local smoothing estimates for the averages $A_t$ over non-degenerate curves to curves lying in the model class. To precisely describe this reduction, it is useful to make the choice of cutoff function explicit in the notation by writing $A_t[\gamma,\chi]$ for the operator $A_t$ as defined in \eqref{eq:averaging operator}.

\begin{proposition}\label{quant nondeg prop} Let $\gamma \colon I \to \R^n$ be a non-degenerate curve, $\chi \in C^{\infty}_c(\R)$ be supported on the interior of $I$ and $0 < \delta \ll 1$. There exists some $\gamma^* \in \mathfrak{G}_n(\delta)$ and $\chi^* \in C^{\infty}_c(\R)$ such that 
\begin{equation*}
    \|A_t[\gamma,\chi]\|_{L^p(\R^n) \to L^p_\alpha(\R^n \times [1,2])} \sim_{\gamma, \chi, \delta, p, \alpha} \|A_t[\gamma^*,\chi^*]\|_{L^p(\R^n) \to L^p_\alpha(\R^n \times [1,2])}
\end{equation*}
for all $1 \leq p < \infty$ and $0 \leq \alpha \leq 1$. Furthermore, $\chi^*$ may be chosen to satisfy $\supp \chi^* \subseteq [-\delta,\delta]$.
\end{proposition}

As a consequence of Proposition~\ref{quant nondeg prop}, it suffices to fix $\delta_0 > 0$ and prove Theorem~\ref{thm:LS nec} in the special case where $\gamma \in \mathfrak{G}_n(\delta_0)$ and $\supp \chi \subseteq I_0 := [-\delta_0,\delta_0]$. Thus, henceforth we work with some fixed $\delta_0$, chosen to satisfy the forthcoming requirements of the proofs. For the sake of concreteness, the choice of $\delta_0 := 10^{-10^5}$ is more than enough for our purposes.




\subsection{The worst decay cone}\label{subsec: cone}

Key to the local smoothing problem is to understand the decay properties of the Fourier transform $\widehat{\mu}_t$ of the underlying measures $\mu_t$ from \eqref{eq: measure}. Here we recap some basic facts in this vein, which have appeared in earlier works such as \cite{PS2007, BGHS_Hel, KLO_Hel, KLO_Sob}.\medskip

By Proposition~\ref{quant nondeg prop} we may assume without loss of generality that $\gamma \in \mathfrak{G}_n(\delta_0)$ for some small $0 < \delta_0 \ll 1$ and that the cutoff $\chi$ in the definition of $A_{t}$ is supported in $I_0 = [-\delta_0, \delta_0]$. Since $\gamma$ is non-degenerate, we have by van der Corput lemma that
\begin{equation}\label{eq:decay}
    |\widehat{\mu}_t(\xi)| \lesssim (1+|\xi|)^{-1/n}
\end{equation}
uniformly in $1 \leq t \leq 2$. In view of the van der Corput lemma, the worst decay cone where \eqref{eq:decay} cannot be improved should correspond to the $\xi$ for which the derivatives $\inn{\gamma^{(j)}(s)}{\xi}$, $1 \leq j \leq n-1$, all simultaneously vanish for some $s \in I_0$. In order to describe this region, first note that 
\begin{equation*}
    \inn{\gamma^{(n-1)}(s_0)}{\xi_0} = 0 \quad \textrm{and} \quad \frac{\partial}{\partial s}\inn{\gamma^{(n-1)}(s)}{\xi} \Big|_{\substack{s=s_0 \\ \xi = \xi_0}} = 1
\end{equation*}
for $(s_0, \xi_0) = (0, \vec{e}_n)$, by the reduction $\gamma^{(j)}(0) = \vec{e}_j$ for $1 \leq j \leq n$. Consequently, provided the support of $\chi$ is chosen sufficiently small, by the implicit function theorem and homogeneity there exists a constant $c > 0$ and a smooth mapping 
\begin{equation}\label{eq: theta}
    \theta \colon \Xi  \to I_0, \qquad \textrm{where} \quad \Xi := \big\{ \xi = (\xi', \xi_n) \in \widehat{\R}^n \backslash \{0\}: |\xi'| \leq c|\xi_n| \big\},
\end{equation}
such that $s = \theta(\xi)$ is the unique solution in $I$ to the equation $\inn{\gamma^{(n-1)}(s)}{\xi} = 0$ whenever $\xi \in \Xi$. Note that  $\theta$ is homogeneous of degree zero.

Further consider the system of $n$ equations in $n+1$ variables given by 
\Be\label{bad cone}
\begin{cases} 
&\inn{\gamma^{(j)}(s)}{\xi}=0 \qquad \text{for $1 \leq j \leq n-1$,}
\\
& \xi_n=1.
\end{cases}
\Ee 
Again, by the reduction $\gamma^{(j)}(0) = \vec{e}_j$ for $1 \leq j \leq n$, this can be solved for suitably localised $\xi$ using the implicit function theorem, expressing $s$, $\xi_1$, ... $\xi_{n-2}$ as  functions of $\xi_{n-1}$. Thus \eqref{bad cone} holds if and only if
\Be\label{solving on the bad cone} 
 \begin{aligned} 
\xi_i&=\Gamma_i(\xi_{n-1}), \qquad 1 \leq i \leq n-2, \\s&=\theta (\Gamma_1(\xi_{n-1}),\dots, \Gamma_{n-2}(\xi_{n-1}), \xi_{n-1},1),
\end{aligned}
\Ee for some smooth functions $\Gamma_i$, $i=1,\dots, n-2$  satisfying $\Gamma_i(0)=0$. 
On $I$ we form the  $\bbR^n$-valued  function $\tau\mapsto \Gamma(\tau)$
with the first $n-2$ components as defined in \eqref{solving on the bad cone} and 
\begin{equation*}
\Gamma_{n-1}(\tau) := \tau, \qquad \Gamma_n(\tau) := 1.
\end{equation*}
With this definition, the formul\ae\, in \eqref{solving on the bad cone} can be succinctly expressed as
\begin{equation*}
    \xi = \Gamma(\xi_{n-1}), \qquad 
    s = \theta\circ \Gamma(\xi_{n-1}).
\end{equation*}
Moreover, the `worst decay cone' can then be defined as the cone generated by the curve $\Gamma$, given by
\begin{equation*}
    \mathcal{C} := \big\{ \lambda \Gamma(\tau) : \lambda > 0 \textrm{ and } \tau \in I \big\}.
\end{equation*}
Our counterexample will live near $\mathcal{C}$ in the frequency domain.




\section{The counterexample}

We now turn to the proof of Theorem~\ref{thm:LS nec}. The necessary condition $\sigma < \min \{\frac{1}{n}, \frac{2}{p}\}$ was already proved in \cite[Proposition 3.9]{KLO_Sob}. We shall only focus on the necessary condition $\sigma < \frac{1}{n}(\frac{1}{2}+\frac{2}{p})$. Given $\lambda > 0$, consider the family of band-limited Schwartz functions
\begin{equation*}
\cZ_\la : =\big\{f\in \cS(\bbR^n): \supp\widehat{f}\subset\{\xi \in \widehat{\R}^n :\la/2\le|\xi|\le 2\la\} \big\}.
\end{equation*}
By elementary Sobolev space theory, it suffices to prove the following proposition.

\begin{proposition}\label{necessity proposition} If $\gamma \colon I \to \R^n$ is a smooth curve satisfying the non-degeneracy hypothesis \eqref{eq:nondegenerate} and $p\ge 2$, then for all $\varepsilon > 0$ we have
\begin{equation*}
   \sup\big\{ \|A_{t}f\|_{L^p(\R^n \times [1,2])}: f\in L^p\cap\cZ_\la, \quad\|f\|_{L^p(\R^n)} = 1\big \} \gtrsim_{p,\gamma, \varepsilon} \la^{-\frac{1}{n}(\frac{1}{2} + \frac{2}{p}) - \varepsilon}. 
\end{equation*}
\end{proposition}

As mentioned in the introduction, the example considered here is a direct modification of that in \cite[Proposition 3.3]{BGHS_Sob}. This amounts to combining a sharp example of Wolff \cite{Wolff2000} for $\ell^p$-decoupling inequalities with a stationary phase analysis of the Fourier multiplier $\widehat{\mu}_{\gamma}$. To prove Proposition~\ref{necessity proposition}, however, we must take into account any smoothing effect from averaging in time; this feature is not present in the analysis in \cite{BGHS_Sob}. The key additional observation is that, for our example, the output function $A_t f$ does not \textit{travel} in a certain time interval of length $\lambda^{-1/n}$.\medskip

Inspired by the example in \cite{Wolff2000}, we consider functions with Fourier support on a union of balls with centres lying on the worst decay cone $\mathcal{C}$. To this end, let $c_0 > 0$ be a small dimensional constant, chosen to satisfy the forthcoming requirements of the argument, and 
\begin{equation*}
   \fN(\la) := \Z \cap \{s \in \R : |s| \le c_0 \la^{1/n}\} .
\end{equation*}
The centres of the aforementioned balls are then given by
\begin{equation*}
\xi^\nu := \la \Gamma  (\nu\la^{-1/n}) \qquad \textrm{for all $\nu \in \fN(\la)$,}
\end{equation*}
 where $\Gamma$ is the parametrisation of the cone $\mathcal{C}$, introduced in \S\ref{subsec: cone}. 

Fix $\eta \in C_c^{\infty}(\widehat{\R}^n)$ satisfying $\eta(\xi) = 1$ if $|\xi| \leq 1/2$ and $\eta(\xi) = 0$ if $|\xi| \geq 1$. 
Let  $0 < \rho < 1$ be another dimensional constant, again chosen small enough to satisfy the forthcoming requirements of the argument, and define Schwartz functions $g_{\nu}$ for $\nu \in \fN(\la)$  via the Fourier transform by
\begin{equation*}
    \widehat{g}_{\nu}(\xi) := \eta\big( \la^{-1/n}\rho^{-1}(\xi - \xi^{\nu})\big).
\end{equation*}
We shall also consider
\begin{equation}\label{eq:g_nu + def}
    \widehat{g}_{\nu,+}(\xi) := \eta_+(\la^{-1/n}\rho^{-1}(\xi-\xi^\nu))
\end{equation}
where $\eta_+\in C^\infty_c(\widehat \bbR^n)$ is such that  $\eta_+(\xi)=1$ for $|\xi|\le 1$ and $\eta_+(\xi)=0$ for $|\xi|>3/2$, so that  $\widehat{g}_{\nu} = \widehat{g}_{\nu,+} \cdot \widehat{g}_{\nu}$.
Further, let 
\begin{equation*}
    \phi(\xi) := \inn{\gamma \circ \theta(\xi)}{\xi} \qquad \textrm{and} \qquad 
    u_n(\xi) :=\inn{\gamma^{(n)} \circ \theta (\xi)}{\xi}.
\end{equation*}
for $\theta(\xi)$ as in \eqref{eq: theta}, and define the constant $\alpha_n$ by
\begin{equation*}
\alpha_{n}:=%
\begin{cases}
\tfrac2n\Gamma(\tfrac1n)\sin(\tfrac{(n-1)\pi}{2n}) \, & \textrm{if $n$ is odd},\\
\tfrac2n\Gamma(\tfrac1n) \exp(i\tfrac{\pi}{2n}) \, & \textrm{if $n$ is even}.
\end{cases}
\end{equation*}

We next record the assymptotics and behaviour of the multiplier $\widehat{\mu}_{\gamma}$ near the support of the $\widehat{g}_\nu$, which is inside the worst decay cone. 

\begin{lemma}\label{asy-lemma} 
Let $1 \leq t \leq 2$. If $c_0$, $\rho > 0$ are chosen sufficiently small, then for all $\la \geq 1$ and $\nu\in \fN(\la)$ the identity
\begin{equation*}
    \widehat{\mu}_{t}(\xi)= e^{-i t\phi(\xi)} m_t(\xi)
\end{equation*} 
holds on $\supp \widehat{g}_{+,\nu}$ where
\begin{enumerate}[i)]
    \item The function $m_t$ satisfies the asymptotics
    \[
    | m_t(\xi) -  \alpha_n \chi(\theta(\xi)) (t u_n(\xi))^{-1/n}| \lesssim_\chi   \rho \lambda^{-1/n} + \lambda^{-2/n} (1+\beta_n \log \lambda)
    \]
     for $\xi \in \supp \widehat{g}_{\nu,+}$; here $\beta_2:=1$ and $\beta_n:=0$ for $n>2$;
    \smallskip
    
    \item The function $m_t$ satisfies the derivative bounds
    \begin{equation}\label{eq:der m}
    |\partial^\alpha_\xi m_t(\xi)| \lesssim \lambda^{-1/n - |\alpha|/n}  \quad \text{ for all $\alpha \in \N_0^n$}
    \end{equation}
    and $\xi \in \supp \widehat{g}_{\nu,+}$.
\end{enumerate}
\end{lemma}

A  variant of this lemma appears in {\cite[Lemma 3.4]{BGHS_Sob}}, where it followed as a consequence of \cite[Lemma 5.1]{BGGIST2007}. One may also deduce Lemma~\ref{asy-lemma}  from \cite[Lemma 5.1]{BGGIST2007}, using similar arguments to those in \cite{BGHS_Sob}; we leave the details to the interested reader.

\begin{proof}[Proof (of Proposition \ref{necessity proposition})] Given $\varepsilon > 0$, it suffices to show the proposition holds for $\lambda$ sufficiently large, depending on $\varepsilon$ and $n$. \medskip 

For each $\nu \in \fN(\la)$ define $f_\nu$ by
 \[ \widehat{f}_\nu (\xi):= \lambda^{1/n}  e^{i \phi(\xi)} \widehat{g}_\nu (\xi) \]
and consider the function
\begin{equation*}
    f :=\sum_{\nu \in \fN(\lambda)} f_\nu.
\end{equation*}
Note that the functions $\widehat{f}_\nu$ are essentially like $\widehat{g}_\nu \cdot (\widehat{\mu}_1)^{-1}$, which were the input functions considered in \cite[p. 11]{BGHS_Sob}. Arguing as in there (in fact, the integration-by-parts is slightly easier here since there is no symbol $a_\nu$) we obtain 
\begin{equation}\label {eq:upper bound f}
\|f  \|_{L^p(\R^n)} \lc \la^{(n+1)/n-(n-1)/np}.
\end{equation}

We will next show that 
\begin{equation}\label{eq:lower bound Af}
    \| A_t f \|_{L^p(\R^n \times [1, 1+\lambda^{-1/n}])} \gtrsim \lambda^{1- \frac{1}{p} + \frac{1}{2n} - \frac{1}{np}-\varepsilon}.
\end{equation}
Assuming this temporarily and combining it with \eqref{eq:upper bound f}, we obtain
\begin{equation*}
    \sup_{f\in L^p\cap\cZ_\la} \frac{\| A_t f \|_{L^p(\R^n \times [1,2])}}{\| f \|_{L^p(\R^n)}} \gtrsim \lambda^{-\frac{1}{p} - \frac{1}{2n} - \frac{1}{np} + \frac{n-1}{np}-\varepsilon} = \lambda^{-\frac{1}{n}(\frac{1}{2}+\frac{2}{p})-\varepsilon}
\end{equation*}
which is the desired bound stated in Proposition~\ref{necessity proposition}.\medskip

Turning to the proof of \eqref{eq:lower bound Af}, we claim that for $|t-1| \lesssim \lambda^{-1/n}$, each piece $A_t f_\nu$ of the operator satisfies the lower bound
\begin{equation}\label{eq:lower bound piece}
    \| A_t f_\nu \|_{L^2(\R^n)} \gtrsim \lambda^{1/2} \quad \text{ for  all $\nu \in \fN(\lambda)$}.
\end{equation}
Once this is proved, Plancherel's theorem (using the disjoint supports of $\widehat{f}_\nu$) and the fact that $\# \fN(\lambda) \sim \lambda^{1/n}$ imply that
\begin{equation}\label{eq: L2 orthogonality}
  \Big\| \sum_{\nu \in \fN(\lambda)} A_t f_\nu \Big\|_{L^2(\R^n)} \gtrsim  \Big( \sum_{\nu \in \fN(\lambda)} \|A_t f_\nu\|_{L^2(\R^n)}^{1/2}  \Big)^{1/2} \gtrsim \lambda^{\frac{1}{2} + \frac{1}{2n}}.
\end{equation}

On the other hand, we also claim that for $|t-1| \lesssim \lambda^{-1/n}$, the function $A_t f$ concentrates in $B(0, \lambda^{-1/n})$. More precisely, for all $\varepsilon > 0$ there exists some $R \geq 1$ such that
\begin{equation}\label{eq:concentration}
    \Big\| \sum_{\nu \in \fN(\lambda)} A_t f_\nu \Big\|_{L^2(\R^n)} \sim \Big\| \sum_{\nu \in \fN(\lambda)} A_t f_\nu \Big\|_{L^2(B(0, \lambda^{-(1-\varepsilon)/n}))}
\end{equation}
holds for all $\lambda \geq R$.\medskip

Once we have verified \eqref{eq:lower bound piece} and \eqref{eq:concentration}, we may apply H\"older's inequality, \eqref{eq:concentration} and \eqref{eq: L2 orthogonality} to deduce that
\begin{align*}
    \Big\| \sum_{\nu \in \fN(\lambda)} A_t f_\nu \Big\|_{L^p(\R^n \times [1,2])}  & \gtrsim \lambda^{\frac{1}{2}-\frac{1}{p} - \varepsilon} \Big( \int_1^{1+\lambda^{-1/n}} \Big\| \sum_{\nu \in \fN(\lambda)} A f_\nu \Big\|_{L^2(B(0,\lambda^{-(1-\varepsilon)/n}))}^p \, \ud t \Big)^{1/p}  \\
    & \sim \lambda^{\frac{1}{2}-\frac{1}{p} - \varepsilon} \Big( \int_1^{1+\lambda^{-1/n}} \Big\| \sum_{\nu \in \fN(\lambda)} A f_\nu \Big\|_{L^2(\R^n)}^p \, \ud t \Big)^{1/p}  \\
    & \gtrsim \lambda^{1-\frac{1}{p} + \frac{1}{2n}-\frac{1}{np} - \varepsilon}
\end{align*}
for $\lambda \geq R$, which is the desired lower bound \eqref{eq:lower bound Af}.\medskip

We turn to the proof of the frequency localised $L^2$ bound \eqref{eq:lower bound piece}. 
Given $\nu \in \fN(\lambda)$ and $ 1 \leq t \leq 2$, let
\begin{equation*}
    c_{t,\nu} := \chi (\theta(\xi^\nu))\frac{\lambda^{1/n}}{(t u_n(\xi^\nu))^{1/n}}
\end{equation*}
and note that $c_{t,\nu} \sim 1$. 
We use Lemma \ref{asy-lemma} i) to write
\begin{equation}\label{eq:A fnu hat dec}
    (A_t f_\nu)\;\widehat{}\;(\xi) =\alpha_n c_{t,\nu} e^{-i(t-1)\phi(\xi)} \widehat{g}_\nu(\xi)+(H_{t,\nu})\;\widehat{}\;(\xi) \widehat{g}_\nu (\xi)
\end{equation}
where
\begin{equation}\label{eq:H def}
    H_{t,\nu}(x) :=  \frac{1}{(2\pi)^n} \int_{\widehat{\R}^n } e^{i \inn{x}{\xi}} e^{-i(t-1)\phi(\xi) } \big( \lambda^{1/n}m_t(\xi) - \alpha_n c_{t,\nu} \big) \widehat{g}_{\nu, +}(\xi)\,\ud \xi
\end{equation}
and $g_{\nu, +}$ is as in \eqref{eq:g_nu + def}.

By the definition of $\theta$ (see \eqref{eq: theta}) and the mean value theorem
\begin{equation*}
    |\chi(\theta(\xi)) - \chi(\theta(\xi^\nu))| \lesssim \lambda^{-1} |\xi-\xi^\nu| \lesssim \rho \lambda^{-(n-1)/n} \quad \text{for $\xi \in \supp \widehat{g}_{\nu, +}$.}
\end{equation*}
Similarly,
\begin{equation*}
    |u_n(\xi)^{-1/n} - u_n(\xi^\nu)^{-1/n}| \lesssim \lambda^{-1-1/n} |\xi - \xi^\nu| \lesssim\rho \lambda^{-1} \quad \text{for $\xi \in \supp \widehat{g}_{\nu, +}$.}
\end{equation*}
Consequently,
\begin{equation*}
    \Big|\chi (\theta(\xi))\frac{\lambda^{1/n}}{( u_n(\xi))^{1/n}}-\chi (\theta(\xi^\nu))\frac{\lambda^{1/n}}{( u_n(\xi^\nu))^{1/n}}\Big| \lesssim \rho \lambda^{-(n-1)/n} \quad \text{for $\xi \in \supp \widehat{g}_{\nu, +}$.}
\end{equation*}
Thus, by Lemma \ref{asy-lemma} i) and the triangle inequality,
\begin{equation}\label{eq:mult approx nu}
    |\lambda^{1/n}m_t(\xi) - \alpha_n c_{t,\nu}| \lesssim_\chi \rho + \lambda^{-1/n} (1+ \beta_n \log \lambda) \quad \text{for $\xi \in \supp \widehat{g}_{\nu, +}$}
\end{equation}
uniformly in $1 \leq t \leq 2$.

Recall also from Lemma \ref{asy-lemma} ii) that
\begin{equation}\label{eq:mt der}
    \lambda^{1/n}|\partial^\alpha_\xi m_t(\xi)| \lesssim_\alpha \lambda^{-|\alpha|/n}= \rho^{|\alpha|} (\rho \lambda^{1/n})^{-|\alpha|} \quad \text{for all $\alpha \in \N_0^n$.}
\end{equation}
Furthermore, by the definition of $\widehat{g}_{\nu, +}$ we trivially have
\begin{equation}\label{eq:g hat der}
    |\partial_\xi^\alpha \widehat{g}_{\nu,+}(\xi)|\lesssim (\rho\lambda^{1/n})^{-|\alpha|} \quad \text{for all $\alpha \in \N_0^n$.}.
\end{equation}
Finally, by the homogeneity of $\phi$ we have that for $|t-1|\lesssim \lambda^{-1/n}$,
\begin{equation}\label{eq:exp der}
    |\partial_\xi^\alpha e^{-i(t-1)\phi(\xi)}| \lesssim |t-1|^{|\alpha|} = \rho^{|\alpha|} (\rho \lambda^{1/n})^{-|\alpha|} \quad \text{for all $\alpha \in \N_0^n$}.
\end{equation}
In view of \eqref{eq:mult approx nu}, \eqref{eq:mt der}, \eqref{eq:g hat der} and \eqref{eq:exp der}, we have by integration-by-parts in \eqref{eq:H def} that
\begin{equation*}
    |H_{t,\nu}(x)| \lesssim_N \rho  \frac{\rho^n \lambda}{(1+\rho \lambda^{1/n} |x|)^N} \quad \text{for $|t-1|\lesssim \lambda^{-1/n}$, $\, N \in \N$.}
\end{equation*}
Thus, recalling \eqref{eq:A fnu hat dec}, there exists a dimensional constant $C > 0$ such that
\begin{equation*}
    \|A_t f_\nu \|_{L^2(\R^n)} \geq  |\alpha_n| c_{t,\nu} \|  g_\nu \|_{L^2(\R^n)} - \|H_{t,\nu} \ast g_\nu \|_{L^2(\R^n)} \geq (|\alpha_n| c_{t,\nu} - C \rho) \| g_\nu \|_{L^2(\R^n)}
\end{equation*}
holds for all $|t-1|\lesssim \lambda^{-1/n}$, where we have used Plancherel's theorem and Young's convolution inequality. Provided $\rho>0$ is chosen sufficiently small,
\begin{equation*}
    \| A_t f_\nu \|_{L^2(\R^n)} \gtrsim \| g_\nu \|_{L^2(\R^n)} \sim \lambda^{1/2}
\end{equation*}
which establishes \eqref{eq:lower bound piece}.\medskip

We now turn to the proof of the concentration estimate \eqref{eq:concentration}. Expressing $\widehat{\mu}_t$ as in Lemma \ref{asy-lemma}, we have
\begin{equation*}
    A_t f_\nu(x)= \frac{1}{(2\pi)^n} \int_{\widehat{\R}^n} e^{i \inn{x}{\xi}} e^{-i(t-1)\phi(\xi)} \lambda^{1/n}m_t(\xi)\widehat{g}_{\nu}(\xi) \, \ud \xi.
\end{equation*}
An integration-by-parts argument, similar to that used above, then shows that
if $|t-1| \lesssim \lambda^{-1/n}$, then 
\begin{equation*}
    |A_t f_\nu (x)| \lesssim_N \frac{\rho^n \lambda}{(1+\rho \lambda^{1/n} |x|)^N} \qquad \textrm{for all $N \in \N$. }
\end{equation*}
Here we have used the definition of $\widehat{g}_\nu$, the derivative bounds \eqref{eq:der m} for $m_t$ in Lemma \ref{asy-lemma} ii), and that the phase function $\phi$ is homogeneous. From this, one readily deduces that
\begin{equation*}
    \Big\| \sum_{\nu \in \fN(\lambda)} A_t f_\nu \Big\|_{L^2(\R^n \setminus B(0, \lambda^{-(1 - \varepsilon)/n}))} \lesssim_{\varepsilon} 1.
\end{equation*}
Combining this with \eqref{eq: L2 orthogonality}, we see that \eqref{eq:concentration} holds for sufficiently large $\lambda$. This concludes the proof of the proposition.
\end{proof}




\bibliography{Reference}

\bibliographystyle{plain}
\end{document}